\documentclass{svproc}
\usepackage{url}
\usepackage{caption}
\usepackage{multirow}
\usepackage{subfigure}
\usepackage{caption}
\usepackage{wrapfig}
\usepackage{epsfig}
\usepackage{subfig}
\usepackage{color}

\def\RR{{\mathop{{\rm I}\kern-.2em{\rm R}}\nolimits}}
\def\I{{\bf I}}
\def\d{{\bf d}}
\def\i{{\bf i}}
\def\p{{\bf p}}
\def\s{{\bf s}}
\def\t{{\bf t}}
\def\X{{\bf X}}
\def\mk{{\mathcal{K}}}
\def\mg{{\mathcal{G}}}
\def\blam{\mbox{\boldmath $\lambda$}}
\def\bmu{\mbox{\boldmath $\mu$}}
\def\cB{{\mathcal{B}}}

\begin{document}
\mainmatter              
\title{Cubature rules based on bivariate spline quasi-interpolation for  weakly singular integrals}
\titlerunning{Spline quasi-interpolation applied to singular cubature}  
%
\author{Antonella Falini\inst{1} \and Tadej Kandu\v{c}\inst{2} \and
Maria Lucia Sampoli\inst{3} \and Alessandra Sestini\inst{4}}
\authorrunning{Antonella Falini et al.} 
%
%
\institute{Dept. of Computer Science, Univ. of Bari, Italy\\
\email{antonella.falini@uniba.it}
\and
 Faculty of  Mathematics and Physics, Univ. of Ljubljana, Slovenia\\
\email{tadej.kanduc@fmf.uni-lj.si}
\and
 Dept. of Information Engineering and Mathematics, Univ. of Siena, Italy\\
\email{marialucia.sampoli@unisi.it}
\and
 Dept. of Mathematics and Computer Science, Univ. of Florence, Italy\\
 \email{alessandra.sestini@unifi.it}}

\maketitle              

\begin{abstract}
In this paper we present a new class of cubature rules with the aim of accurately integrating weakly singular double integrals. In particular we focus on those integrals coming from the discretization of Boundary Integral Equations for 3D Laplace boundary value problems, using a collocation method within the Isogeometric Analysis paradigm. In such setting the regular part of the integrand can be defined as the product of a tensor product B-spline and a general function. The rules are derived by using first the spline quasi-interpolation  approach to approximate  such function and then  the extension of a well known algorithm for spline product to the bivariate setting.  In this way efficiency is ensured, since  the locality of any spline quasi-interpolation scheme is combined with the capability of an ad--hoc treatment of the B-spline factor. The numerical integration is performed on the whole support of the B-spline factor by exploiting inter-element continuity of the integrands.

\keywords{Cubature rules, Singular and nearly singular integrals, Boundary Element Methods, Tensor product B-splines, Spline quasi-interpolation, Spline product, Isogeometric Analysis.}
\end{abstract}
\section{Introduction}
The accurate and efficient numerical evaluation of singular integrals is one of the crucial steps in the numerical simulation of differential problems that can be modeled by Boundary Integral Equations (BIEs)\cite{HsiaoWendl08}. This is the case when relying on Boundary Element Methods (BEMs), which were introduced in the eighties for the numerical solution of several differential problems, either stationary and evolutive, see for example \cite{Wend85,Costabel94} and references therein. The main features of BEMs are the reduction of the problem dimension  and the easiness of application to problems on unbounded domains.
On the other hand it is well known that one of the major efforts with any BEM formulation consists in having to deal with singular and nearly singular integrals, which require special numerical treatment in order to preserve the theoretical convergence order of the numerical solution produced by the adopted discretization.

In this paper we focus on cubature rules for weakly singular integrals. Since the interest in integrals of this kind comes from the isogeometric formulation of BEMs,  let us   briefly recall  their main ideas. 
The first formulation of BEMs considered a piecewise linear approximation of the boundary of the domain, but more accurate curvilinear BEMs already appeared in the nineties. In the latter methods the boundary of a 2D domain is described through  a planar parametric curve. In the parameter domain  of the curve a set of Lagrangian functions  is defined  for the discretization of the considered BIE. The basis of the discretization space where the missing Cauchy data are approximated is just obtained by lifting such functions to the physical boundary of the domain using its parametric representation. Such methodology  is common to collocation and Galerkin approaches and can be extended also to the  isogeometric formulation of a BEM. This is characterized by the significant assumption that the boundary is parametrically represented in B-spline or NURBS form and  the discretization  space $V$ is defined  through B-splines  instead of Lagrangian functions.  This makes possible to increase  the smoothness of functions belonging to $V$ at desired joints between adjacent elements, often guaranteeing a remarkable reduction of the number of  degrees of freedom necessary to attain a certain level of accuracy \cite{ADSS16}. Note that additional flexibility can be achieved by relying on generalized B-splines, see for example \cite{MPS11} and references therein, that can be used for the description of  the geometry and/or  the definition of the discretization space $V$ \cite{ADSS17}.
Furthermore, it has been already shown in the literature that for a 2D IgA--BEM the element--by--element assembly strategy is not anymore strictly necessary  \cite{CMAME18}. This  computational advantage is obtained since the required integrals, even when singular, can be approximated by rules  formulated directly on the support of the B-spline explicitly appearing in the integrand as one of the basis functions generating $V$ \cite{JCAM18}.

The literature on numerical approximation of singular integrals is quite vast and it is difficult to cover all the results on this issue, see for instance the book \cite{Sladek2'98} or the more recent paper \cite{Gao2010} and references therein. As our interest for singular integrals directly descends from their occurrence within the Isogeometric formulation of BEMs (IgA--BEMs), we limit our attention to the integrals of this kind arising in 3D problems.
Singularity removal is often proposed for the numerical treatment of the occurring multivariate weakly singular integrals.  For example in \cite{KFK09} where the 3D Stokes problem is considered,  the singularity is removed by exploiting carefully chosen known solutions of the analyzed partial differential equation. In other papers these integrals are reformulated by using a suitable coordinate transformation, see for example \cite{Tan_etal2019} for Duffy and \cite{Hugh16} for polar transformations. In these cases the additional emerging transformation term approximately cancels out the singularity of the kernel and the resulting integrals become regular. In \cite{GongDong17} an adaptive Gaussian quadrature rule is presented and it is shown that it is able to tackle singular and also near singular integrals. 
However all these approaches do not exploit the smoothness of B-splines, taking only into account their piecewise polynomial nature. For this reason, the related cubature rules are always applied after splitting the integration domain into elements with a consequent increase of the computational cost. Instead, in this paper, the B-spline factor is explicitly treated and the cubature rule is applied on the whole B-spline support, not suffering from inter-element smoothness decrease of B--splines. The rules here proposed are an extension to the bivariate setting of the quadrature formulas for singular integrals introduced in  \cite{JCAM18}. Their key ingredients  are  a  spline quasi-interpolation approach and  the spline product formula  \cite{Morken91}, both considered in their tensor--product formulation.  By exploiting the integration on the whole B-spline support, they are attractive for IgA-BEM also in the 3D case, where a replacement of element-by-element assembly with a function-by-function strategy is even more advantageous.

The paper is organized as follows. First we introduce cubature rules for weakly singular integrals,  showing their effectiveness when the considered kernel is multiplied by a general function and a B-spline.
Then the combination with suitable multiplicative or subtractive techniques specific of the 3D setting  is analyzed, in order  to show that they become applicable to deal with specific singular integrals of interest in the IgA-BEM setting.

\section{The problem}
In this paper we focus on cubature rules   for singular integrals of the following type,
\begin{equation} \label{int}
\int_{R_\I} \mk(\s\,,\,\t) \, B_{\I,\d}(\t)\, f_\s (\t) \, d \t \,, \qquad \s \in R_{\I}^E\,,
\end{equation}
where  $B_{\I,\d}$ is an assigned bivariate B-spline of bi--degree $\d:=(d_1,d_2)$ with support in the rectangle $R_\I$, $R_\I^E \supset R_\I$, and
\begin{equation} \label{kernel}
  \mk(\s\,,\,\t)  :=  \frac{1}{\sqrt{(\t-\s)^T A(\s) (\t-\s)}},\quad \t = (t_1,t_2)\,, \quad \s = (s_1,s_2)\,,
  \end{equation}
with $A(\s)$ denoting a symmetric and positive definite matrix (which ensures that the singularity appears just at $\t = \s$).
Concerning the smoothness requirements for $f_s,$ since our rules are based on the tensor product formulation of (a variant of) an Hermite quasi-interpolation scheme, it is reasonable to assume  $f_s$ belonging to $C^{1,1}(R_\I),$ that is to the space of bivariate functions $g$ such that $\frac{\partial^{i+j} g}{ \partial t_1^i \partial t_2^j}$ is continuous in $R_\I$ for $i,j \le 1.$ 
We refer to \cite{Schumaker} for an introduction on basic properties and definitions of B-splines and in particular on their tensor product bivariate extension. We observe that for $\s \in R_\I$ the  integral in (\ref{int}) is weakly singular and it becomes nearly singular when $\s \in R_\I^E \setminus R_\I,$ with the maximal distance from $R_\I$ of  $ \s \in R_\I^E~\setminus~R_\I$ sufficiently small  to exclude regular integrals. This is in contrast to other approaches proposed in the literature (see for instance \cite{Scuderi09}), where typically different integration methods are used for singular and nearly singular integrals. We also note that our rules numerically compute the integral in (\ref{int}) by approximating only the factor $f_\s$. This is particularly useful when the function $f_\s$  is more regular in $R_I$ than $B_{\I,\d}$, since usually it can be better approximated than the whole product $B_{\I,\d} f_\s$  \cite{JCAM18}.

We outline that the kernel $\mk$ is of interest for BEMs   when $A(\s)$ is the matrix containing the coefficients at $\t = \s$ of the first fundamental form associated to a differentiable parametric surface $\X = \X(\t), \t \in {\cal D} \subset  \RR^2,$
\begin{equation} \label{Adef} A(\t) = \left[\begin{array}{cc} (\X_{t_1} \cdot  \X_{t_1} ) (\t) &  (\X_{t_1} \cdot  \X_{t_2} ) (\t)\cr
 (\X_{t_1} \cdot  \X_{t_2} ) (\t)&  (\X_{t_2} \cdot  \X_{t_2}  ) (\t) \cr\end{array} \right]\,.\end{equation}
Indeed in this case  the quadratic homogeneous polynomial
 \begin{equation} \label{quad} P_\s(\t) := (\t-\s)^T A(\s) (\t-\s)
 \end{equation} collects the  lowest order non-zero terms of the Taylor expansion at $\t = \s$ of $\Vert \X(\t) - \X(\s) \Vert_2^2.$ So $\mk(\s,\t)$ is a local approximation of
\begin{equation} \label{kernelG}
  \mg(\s\,,\,\t)  :=  \frac{1}{ \Vert \X(\t) - \X(\s) \Vert_2 }\,,
  \end{equation}
which is, up to a multiplicative constant, the kernel appearing in the single layer potential,
\begin{equation} \label{original}
 \int_{R_\I} \mg(\s\,,\,\t) \, B_{\I,\d}(\t)\, g_\s(\t)\, d \t\,,
 \end{equation}
for 3D Laplace problems, written in intrinsic coordinates. The B-spline factor in (\ref{original}) corresponds to a basis function of the tensor product spline space V used for the discretization, while  $g_{\s}$ appears in the formulation as the Jacobian of the domain transformation to the parametric domain.
Note that   $ \mg$ is substantially the kernel associated also with the Helmholtz problem, missing only an additional regular trigonometric factor appearing in the fundamental solution of such equation.

In this work  we consider the so-called singularity extraction procedure, based on either a subtractive or a multiplicative technique, to derive a more convenient formulation of the singular integral.
Following this procedure, the integral in (\ref{original}) is transformed into an integral with the same kind of singularity but with a more standard kernel, possibly added to a regular integral.

Denoting with $\mg_a$  the approximating kernel having the same kind of singularity of $\mg$ at $\t = \s,$  with the subtractive technique the integral in (\ref{original}) is decomposed in the following sum,
\begin{equation} \label{sum}
  \int_{R_\I} \mg_a(\s\,,\,\t) \, B_{\I,\d}(\t)\, g_\s(\t)\, d \t + \int_{R_\I} \left(\mg(\s\,,\,\t) - \mg_a(\s\,,\,\t)\right) \, B_{\I,\d}(\t)\, g_\s(\t)\, d \t
\end{equation}
where the second integral is regular if $\mg_a$ is suitably defined. The first integral in (\ref{sum}) is still weakly singular and it becomes equal to the integral in (\ref{int}) if   $\mg_a = \mk$ is chosen and $f_\s = g_\s$ is set. In this case the regularity of $f_\s$  is that of the Jacobian of $\X.$ Then, considering the IgA paradigm, we can observe that it can be low (anyway at least $C^{1,1}$ if $\X$ is a regular $C^{2,2}$ NURBS parameterization) only at the original knots involved in the CAGD representation of $\X$, and not at the other knots used to define the discretization space $V$. Furthermore, without loss of generality, we can assume that the original knots have maximal multiplicity, so that the possible reduction of regularity of $f_\s$ can appear only at the boundary of $R_I.$
With the multiplicative technique, setting $\rho_\s(\t) := \mg(\s\,,\,\t)/\mg_a(\s\,,\,\t),$ and $f_\s(\t) := \rho_\s(\t)\,g_\s(\t),$ we obtain
\begin{equation} \label{mult}
\int_{R_\I} \mg(\s\,,\,\t) \, B_{\I,\d}(\t)\, g_\s(\t)\, d \t  = \int_{R_\I} \mg_a(\s\,,\,\t) \, B_{\I,\d}(\t)\, f_\s(\t) \, d \t\,,
\end{equation}
where the function $f_\s$ is regular, again if $\mg_a$ is suitably defined.
If in particular $\mg_a = \mk,$  we get
\begin{equation} \label{rhodef}
\rho_{\s}(\t)  =   \frac{\sqrt{(\t-\s)^T A(\s) (\t-\s)}}{ \Vert \X(\t) - \X(\s) \Vert_2 }\,,
\end{equation}
 with $A$ defined as in (\ref{Adef}). Note that this reformulation of the singular integral in (\ref{original})  can be considered as a bivariate generalization of the standard one proposed in the literature for dealing with univariate singular kernels, where $\mg_a$ is just defined as  $\mg_a(s,t) = 1/|s-t|$.
In the bivariate setting the function $\rho_\s$ defined in (\ref{rhodef}) is continuous at $\t =\s$, since it can be verified that  $\lim_{\t \rightarrow \s}\rho_\s(\t)$ exists and is equal to  $1.$ Unfortunately $\rho_\s$ is not smoother than $C^0$ at such point for a general surface $\X.$  Thus, when the integral of interest is that defined in (\ref{original}) and $\X$ is a general surface, we would need to consider higher order approximations of $\mg$ instead of $\mk,$ in order to deal with functions $f_\s$ more regular at $\t = \s$ when they are obtained by using the multiplicative technique.  Note that also adopting the subtractive technique this can be useful to increase the regularity of the integrand of the regular integral in (\ref{sum}). To keep the presentation of our rules concise, this technical but important aspect is not addressed in this paper.

\section{Cubature rules based on tensor-product spline quasi-interpolation}
Quasi-Interpolation (QI) is a general approach for approximating a function or a given set of discrete data with low computational cost, see for instance   \cite{Sablo05} and references therein. For a chosen finite dimensional approximating space and a suitable local basis generating it, the coefficients of the approximation are locally computed with explicit formulas by using linear functionals depending on the function and possibly also on its derivatives and/or integrals.
Since there is already an explicit B-spline factor in the considered integral in (\ref{int}), it is particularly beneficial for us to approximate the function $f_s$ using a spline quasi-interpolation operator. That way the B-spline factor is preserved in the expression for the numerical integration and the spline product algorithm can be readily applied \cite{Morken91}. \\
The easiest extension of a univariate QI scheme to the bivariate setting relies on its tensor-product formulation which anyway performs function approximation on a rectangular domain, requiring information at the vertices of a quadrilateral grid of the domain. We add that in the bivariate spline setting  there has recently been  a lot of interest for QI schemes on special type triangulations or  even on general ones adopting macroelements, see for example \cite{BDIR19,GSCAGD18} and references therein. However, since for application to cubature   the analytic expression of the function to be approximated is available and our integration domain is rectangular, for our purposes the tensor-product extension  is more suitable. In particular we adopt a  tensor-product derivative free QI scheme which is a natural choice for numerical integration.

Denoting with $S_{p.T}$ the space of univariate splines with degree $p$ and with $T$ the associated extended knot vector defined  in the reference domain $[-1\,,\,1],\,$ -- $ T = \{ \xi_0\le \cdots,\xi_{p-1}\le \xi_p \le \cdots \le \xi_{m+1} \le \cdots \le \xi_{m+p+1}\}\,,$ with $\xi_j < \xi_{j+p+1}$ and $\xi_p = -1, \xi_{m+1} = 1$ --  a spline $\sigma \in S_{p,T}$ can be represented by using the standard B-spline basis, $\cB_{j,p}, j=0,\ldots,m,$
$$\sigma(\cdot) = \sum_{j=0}^m \lambda_j \ \cB_{j,p}(\cdot)\,.$$
Thus a univariate derivative free QI scheme to approximate a univariate function $f$ can be compactly written as follows,
\begin{equation} \label {univQI}
  \blam = C {\bf f}\,,
  \end{equation}
where $\blam := (\lambda_0,\ldots,\lambda_m)^T$ is the vector of the spline coefficients; $C$ is a $(m+1) \times(K+1)$ banded matrix characterizing the scheme; ${\bf f} := (f(\tau_0,)\ldots,f(\tau_K))^T \,,$ with $-1 \le \tau_0 < \cdots < \tau_K \le 1$ completing the characterization of the scheme. On this concern observe that, if $C_{i,j}\,  j=i-L,\ldots,i+U$ are the non vanishing elements in $C,$ it must be required that $\tau_{i-L},\ldots,\tau_{i+U}$ belong to the support of $\cB_{i,p}.$ Furthermore a certain polynomial reproduction capability of the scheme must be required to ensure a suitable convergence order.

\noindent
Within this kind of QI schemes, we refer to the derivative free variant of the Hermite QI method   introduced  in \cite{MS09}. Such variant  requires in input only the values of $f$ at the spline breakpoints, since the derivative values required in the original scheme are approximated with suitable finite differences \cite{MS09}.

In the tensor product formulation of the scheme we have to define a spline  $\sigma $ in the space $ S_{p_1,T_1} \times S_{p_2,T_2},$
$$ \sigma(t_1,t_2) =  \sum_{i=0}^{m_1}  \sum_{j=0}^{m_2} \lambda_{i,j} \  \cB_{i,p_1}(t_1) \cB_{j,p_2}(t_2)\,.$$
Setting $\t := (t_1,t_2)$ and $\I  := \{ (i,j), i=0,\ldots, m_1, j=0,\ldots,m_2\}$ we can compactly write
$$\sigma(\t) = \sum_{\i \in \I} \lambda_\i \, {\cal B}_{\I,\p}(\t)\,,$$
where  ${\cal B}_{\I,\p}(\t) := {\cal B}_{i,p_1}(t_1)  {\cal B}_{j,p_2}(t_2).$ Using for example the lexicographical ordering for the elements of $\I$ and the  Kronecker product between matrices, the tensor product extension of the scheme can be expressed as follows,
\begin{equation} \label {bivQI}
  \blam = (A_1 \otimes A_2) {\bf f}\,,
  \end{equation}
  where  now ${\bf f} = \left(f(\tau_0^{(1)},\tau_0^{(2)})\,,\,f(\tau_0^{(1)},\tau_1^{(2)}),\cdots\,,\,
  f(\tau_{K_1}^{(1)},\tau_{K_2}^{(2)}) \right)^T$ with $f$ denoting a bivariate function   and  $\blam$ is the vector $\blam :=\left( \lambda_{(0,0)},\, \lambda_{(0,1)},\ldots,\lambda_{(m_1,m_2)} \right)^T.$

 In order to extend to the bivariate setting the quadrature rule for singular integrals containing a  B-spline weight developed in \cite{JCAM18}, we need two additional ingredients: a bivariate generalization of the spline product formula and explicit analytical formulas to compute specific singular integrals.
 In more detail, we first consider the tensor product generalization of the algorithm in \cite{Morken91} to express the product $\sigma \ B_{\I,\d}$ in the bivariate B-spline   basis of the product space. Such space has bi--degree $(p_1+d_1,p_2+d_2)$ and the related extended knot vectors in each coordinate direction are obtained by merging $T_k$ and ${\cal T}_k,$ for $k=1,2\,,$ knot vectors in each direction $k$ for $B_{\I}$ and $\sigma,$ respectively.
The other necessary step for approximating the integral in (\ref{int}) consists in the computation of the so-called {\it modified moments},
$$\mu_\i(\s) := \int_{R_\I} {\cal K}(\s\,,\,\t)\, B^{(\Pi)}_\i(\t)\, d\t\,, \qquad \i \in \I^{(\Pi)}\,,$$
where $B^{(\Pi)}_\i, \i \in \I^{(\Pi)},$ denotes the B--spline basis of the product space.
For this aim we need again to generalize to the bivariate setting the univariate recursion for B-splines whose usage in this context was introduced in \cite{CMAME18}. We refer to \cite{BEM3D} for more details on these two steps.

The final approximation of the integral in (\ref{int}) is then simply given by the product $\bmu(\s)^T \blam^{(\Pi)}, $ where $\bmu(\s)$ is the vector containing the above modified moments ordered in lexicographical way and $ \blam^{(\Pi)}$ is a vector of the same length whose entries are  the coefficients expressing $\sigma B_{I,\d}$ in the B-spline basis of the product space.

\section{Numerical Results}
This section is devoted to check the performance of our cubature rules.\\
In the experiments we always assume that the bi-degree $\d =(d,d)$ of the B-spline factor in the integrand of (\ref{int}) is  equal to $(2,2)$ or $(3,3)$ and that $R_{I}=[-1,\ 1]^2$. For simplicity, we consider  a uniform  distribution of the $d+1$ breakpoints of the B-spline in each coordinate direction. In order to deal either with nearly singular and singular integrals, we consider the source points $\s=(s_1,s_2) \in {\cal S}^2 $ with ${\cal S} := \{-1.1,-1,-0.5,0,0.5,1,1.1\}$.\\
The tests are performed on a uniform $N \times N$ grid for the breakpoints of the quasi-interpolating spline $\sigma$, with $N$ ranging from $6$ to $14$ with step $2.$ The bi-degree $\p =(p,p)$ of the quasi-interpolant is set to $(2,2)$ or $(3,3)$.\\

\noindent{\it Example 1}\\
In the first example we consider the quadratic bivariate polynomial function $f_{\s}(\t)=f(\t)=t_1^2+t_2^2$.  The aim of the test is to check the exactness of the proposed cubature rule, since the integration rule is based on the chosen tensor product QI scheme, which is exact on polynomials of bi-degree $(\ell_1,\ell_2)$ with $\ell_k \le p$. For this example the matrix  $A$ defining the kernel ${\cal K}$ in (\ref{kernel}) is just a constant matrix with all unit entries. We verified that already with $N=6$ we get a maximum relative error of $1.54\mbox{e-}13$ for   $\s \in {\cal S}^2$ restricted to the interior of $R_I$. It becomes $7.56\mbox{e-}12$, and $9.60\mbox{e-}12$ when $\s \in {\cal S}^2$ is restricted to the  boundary of $R_I$ and to  values external to $R_I,$ respectively.\\

\noindent{\it Example 2}\\
In order to check the convergence order, in this example we consider $A $ equal to the identity and the analytic function  $f_{\s}(\t)=f(\t)=\exp(t_1t_2)$. The results are collected in Table \ref{exp_table}, where in particular the maximal absolute  errors {\tt errmax1, errmax2} and {\tt errmax3}  are reported, varying the number $N \times N$ of cubature nodes uniformly distributed in $R_I.$ The results show a very good behavior of the rules for the considered test function and matrix.

\begin{table} [htp]
\centering
\resizebox{\columnwidth}{!}{%
\begin{tabular} {|c || c c| c c| c c || c c| c c| c c|}
\hline
\multicolumn{1}{| c |}{$d=2$} &  \multicolumn{6}{||c||}{ {p=2}}&\multicolumn{6}{|c|}{{p=3}} \\
\hline
{$N$}& {\tt errmax1}& {${o}_1$} &{\tt errmax2}& {${o}_2$}&
{\tt errmax3}& {${o}_3$}& {\tt errmax1}& {${o}_1$} &{\tt errmax2}&
{${o}_2$}& {\tt errmax3}& {${o}_3$}\\
\hline
$6$ & 2.5704e-05 &  -- & 4.3428e-05  & -- & 8.3210e-05 &  -- &   1.0520e-06&   -- &2.1322e-06 & -- &2.1322e-06 & --\\
$8$ & 8.4609e-06 &  3.9 & 1.6115e-05  & 3.5 & 1.6697e-05 &  5.6 &   2.7380e-07&  4.7 & 5.4119e-07 & 4.8 &5.4278e-07 & 4.8\\
$10$& 3.6045e-06&  3.8 & 6.9256e-06  & 3.8 & 6.9256e-06 & 3.9 &   9.9469e-08&  4.5 &1.9417e-07 & 4.6 &1.9417e-07 & 4.6\\
$12$& 1.7283e-06 &  4.0 & 3.3031e-06  & 4.1 & 3.3031e-06 & 4.1 &   4.4251e-08&  4.4 &8.5289e-08 & 4.5 &8.5289e-08 & 4.5\\
$14$& 9.1746e-07 & 4.1 & 1.7456e-06 &4.1 & 1.7456e-06& 4.1&
2.2321e-08& 4.4 & 4.2435e-08& 4.5& 4.2435e-08 & 4.5\\
\hline
\hline
\multicolumn{1}{|c|}{$d=3$} &  \multicolumn{6}{||c||}{ {p=2}}&\multicolumn{6}{|c|}{{p=3}} \\
\hline
$6$ &5.0578e-06 &  -- & 1.5198e-05  & -- & 2.5845e-05 &  -- &   3.3475e-07&   -- &8.3595e-07 & -- &8.3595e-07 & --\\
$8$ & 2.6660e-06 &  2.2 & 5.9122e-06  & 3.3 & 5.9122e-06 &  5.1 &   8.7285e-08&  4.7 & 2.1109e-07 & 4.8 &2.1156e-07 & 4.8\\
$10$& 1.1965e-06&  3.6 & 2.6836e-06  & 3.5 & 2.6836e-06 & 3.5 &   3.1949e-08&  4.5 &7.6082e-08 & 4.6 &7.6082e-08 & 4.6\\
$12$& 5.7522e-07 &  4.0 & 1.2883e-06  & 4.0 & 1.2883e-06 & 4.0 &   1.4385e-08&  4.4 &3.3872e-08 & 4.4 &3.3873e-08 & 4.4\\
$14$& 3.0410e-07 & 4.1 & 6.8169e-07 &4.1 & 6.8170e-07& 4.1&
1.0270e-08& 2.2 & 1.7292e-08& 4.4 & 1.7292e-08 & 4.4\\
\hline
\end{tabular}%
}
\vskip 11pt
\caption{Example 2. Maximal absolute cubature error and convergence order for  $\s \in {\cal S}^2$ outside ({\tt errmax1}, ${o}_1$), on the boundary ({\tt errmax2}, ${o}_2$) and inside ({\tt errmax3}, ${o}_3$) the integration domain $R_I,$ for $p=2,3$ and $d=2,3$. } \label{exp_table}
\end{table}

\noindent{\it Example 3}\\
This example considers the case of  the matrix $A$ defined as in (\ref{Adef}), with $\X$ being the standard parameterization for   the lateral surface of a cylinder of radius $r = 2$
$$\X(\t) = \left( r\cos(\pi t_1/4)\,,\,r\sin(\pi t_1/4)\,,\,t_2\right)\,,$$
which implies that   $R_I$  is mapped to a quarter of the lateral cylindrical surface with height $2.$
The factor $f_\s$ in (\ref{int}) is assigned as the product between $\rho_{\s}$  which is  defined in (\ref{rhodef}) and the Jacobian $J(\t),$ with
\begin{equation} \label{Jdef}
J(\t) := \Vert \X_{t_1}(\t) \times \X_{t_2}(\t) \Vert_2.
\end{equation}
This means that the integral with the form in (\ref{int}) considered for this experiment has been obtained from (\ref{original}) by using the multiplicative strategy introduced in (\ref{mult}) with ${\cal G}_a = {\cal K},$ obtaining in this case a $C^{1,1}$ smooth function $\rho_\s$   also when $\s \in R_I.$

\begin{figure}[ht]
\centering
\subfigure[]{\includegraphics[width=6cm]{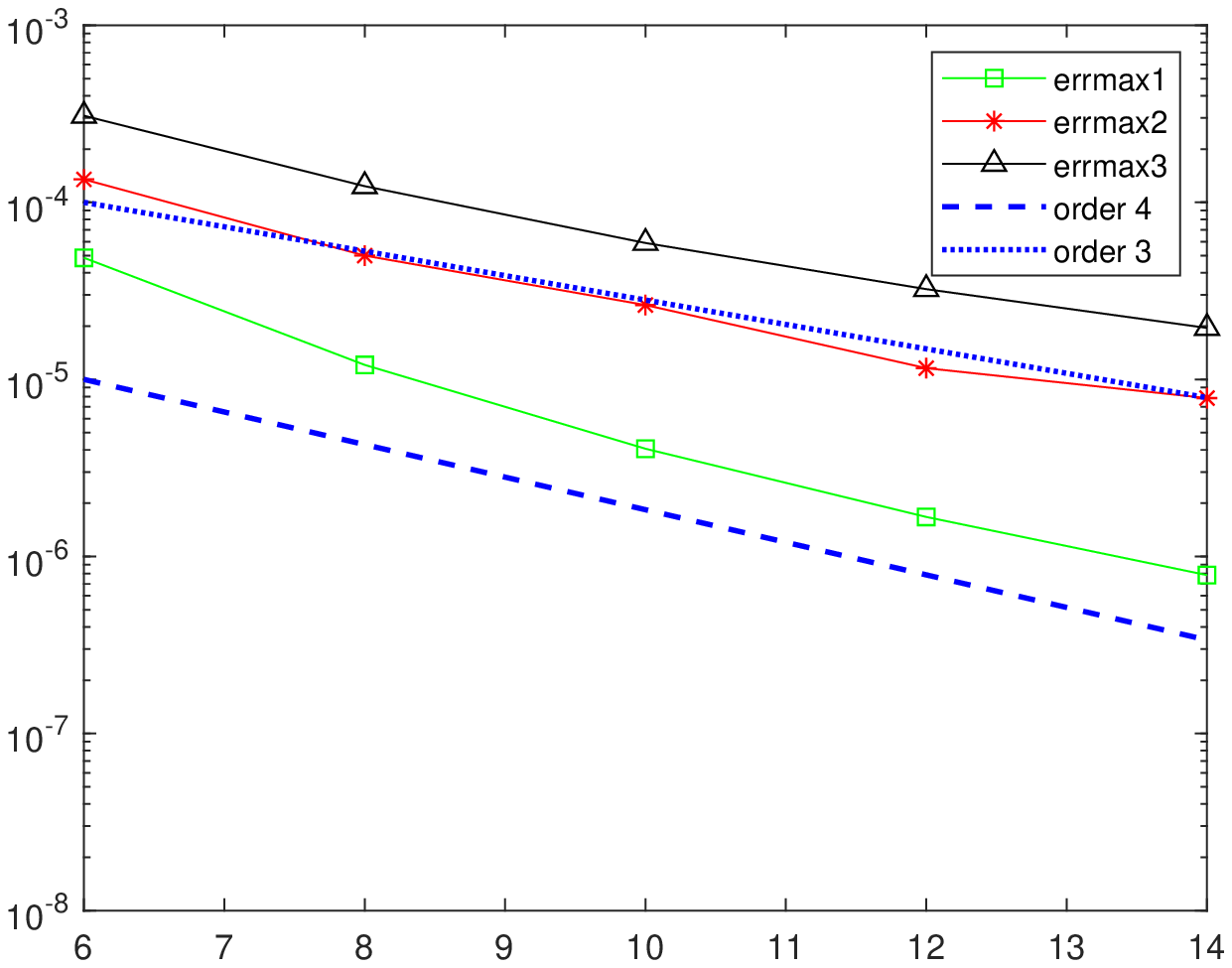}}
\subfigure[]{\includegraphics[width=6cm]{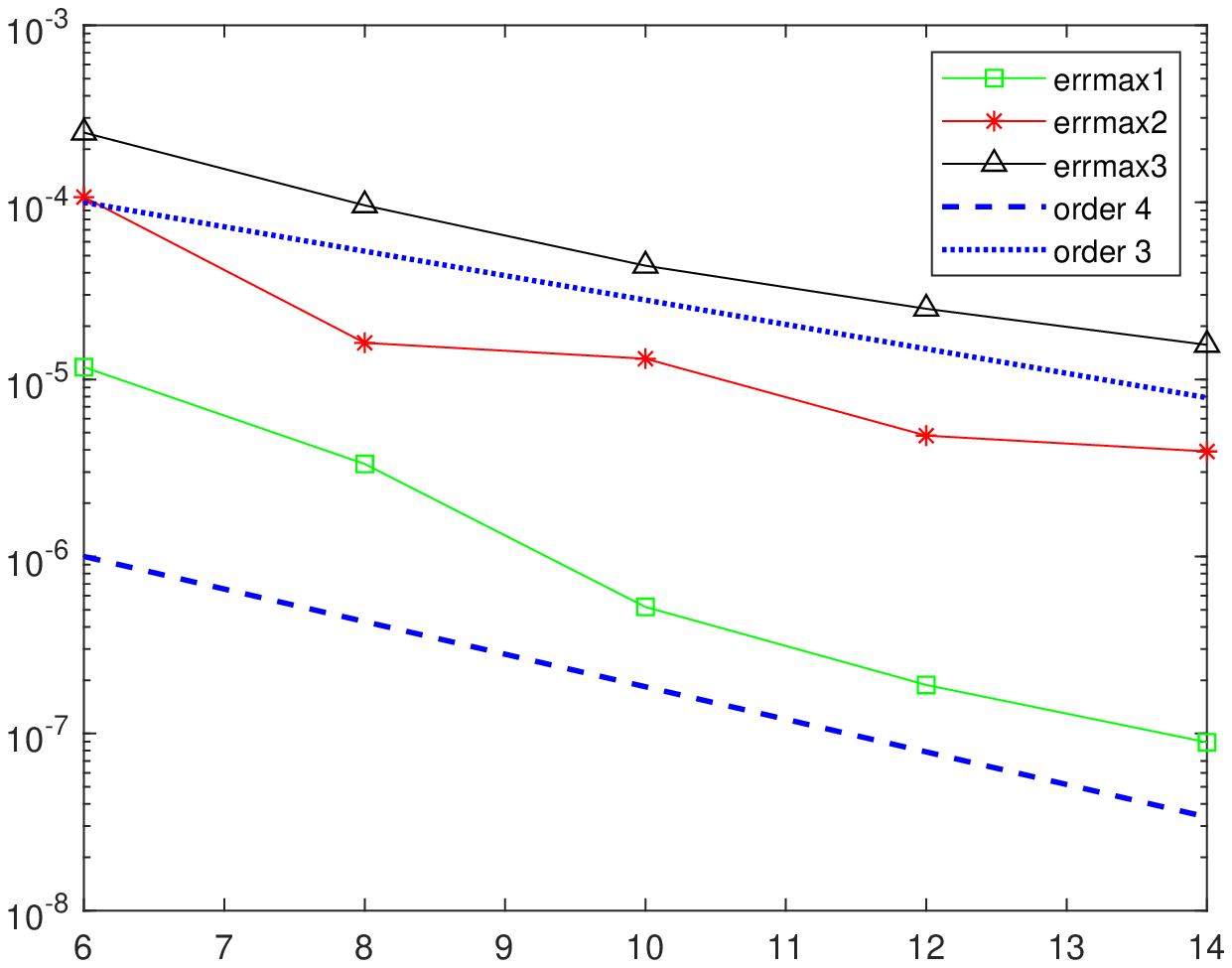}}
\subfigure[]{\includegraphics[width=6cm]{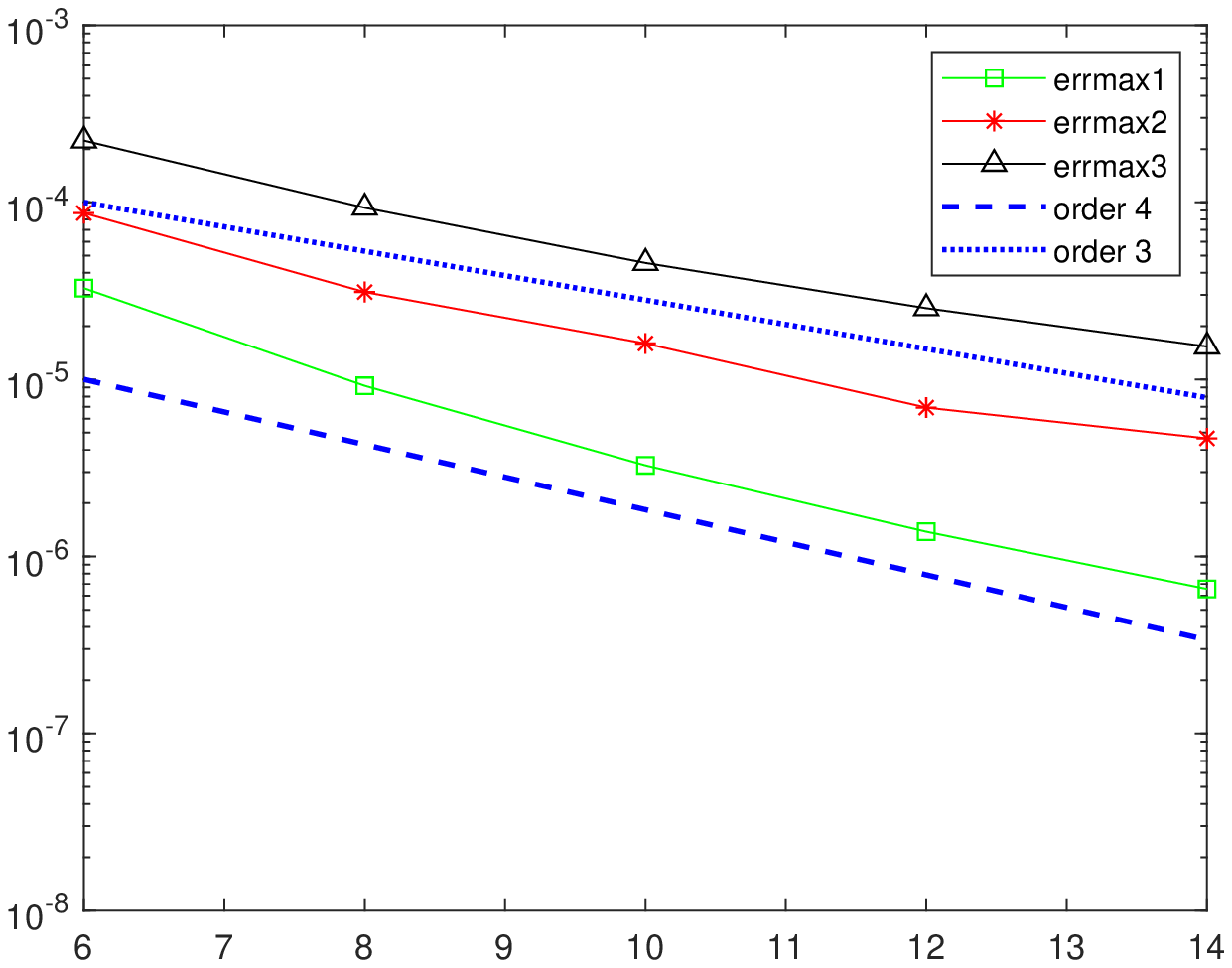}}
\subfigure[]{\includegraphics[width=6cm]{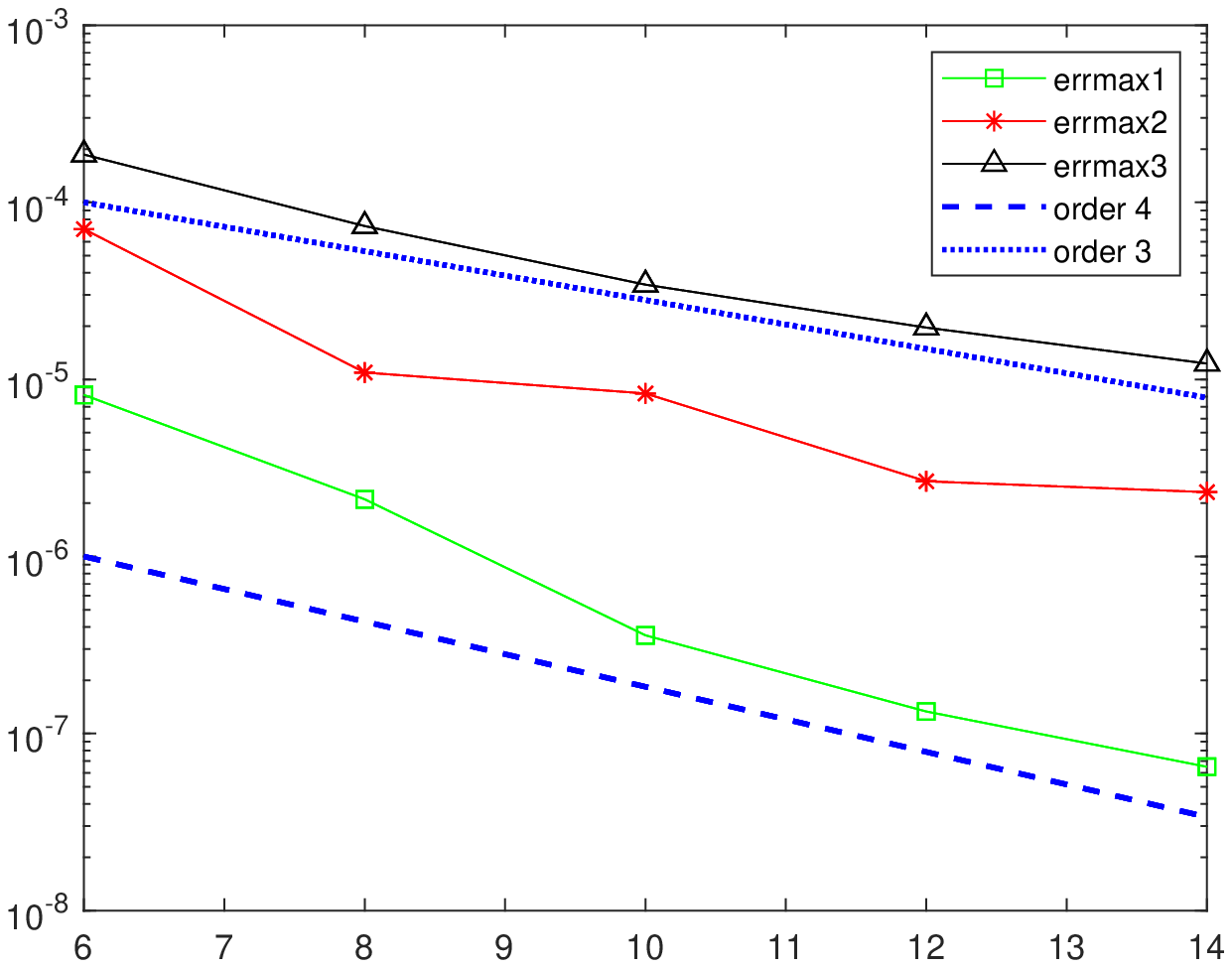}}

\caption{Example 3. The convergence behavior of the absolute cubature errors {\tt errmax1}, {\tt errmax2} and {\tt errmax3} for $d = 2, \ p = 2,$ (a),  $d = 2 , \ p = 3$ (b),
 $d = 3, \ p = 2$ (c) and  $d = 3, \ p = 3$ (d).}
\label{figEx3}
\end{figure}

Figure \ref{figEx3} shows the convergence behavior of the absolute cubature errors {\tt errmax1}, {\tt errmax2} and {\tt errmax3} for the four considered choices of the pair $(d,p).$ Comparing left and right images of the figure and first  referring to {\tt ermax2} and {\tt errmax3} (i.e. when the rules are applied to singular integrals), we can observe that there is not significant advantage in using $p=3$ instead of $p=2$, either from the point of view of the convergence order or from that of the initial ($N=6$) and final ($N=14$)  accuracy.  This is a different behavior with respect to Example 2 where the function $f$ was highly smooth everywhere. Referring to {\tt errmax1} (i.e. for nearly--singular integrals) however, this comment does not hold anymore.  \\
We observe that for the maximum considered value of $N,$ $N=14,$ we achieve a value for {\tt ermax3} of the order of $10^{-5}$ which corresponds to a relative error of the same order; at a first sight this could seem not satisfactory but we remark that the portion of the cylindrical surface taken into account for the integration is quite large. Indeed, repeating the experiment mapping $R_I$ to a smaller portion of the surface, the relative error decreases.
Finally, comparing top and bottom images we can also conclude that  different regularity of the B-spline factor in (\ref{int}) associated with different choices of $d$ does not significantly influence the accuracy of our rules.\\

\noindent{\it Example 4}\\
In the last example, we consider an integral of interest for the BIE formulation of the 3D Helmholtz problem $\Delta u + k^2 u = 0,$ where $k$ is the wave number defined as $k = 2\pi/\lambda,$ with $\lambda$ denoting the wavelength of the electromagnetic radiation. The boundary of the domain of the differential problem is assumed equal to a section of a one sheet hyperboloid which can be parametrically represented as follows,
$$ \X(\t) = (\cos(\pi t_1/4) \sqrt{1+t_2^2}\,,\,\sin(\pi t_1/4 ) \sqrt{1+t_2^2}\,,\,t_2)\,.$$
\begin{figure}[bht]
\centering
\subfigure[]{\includegraphics[width=6cm]{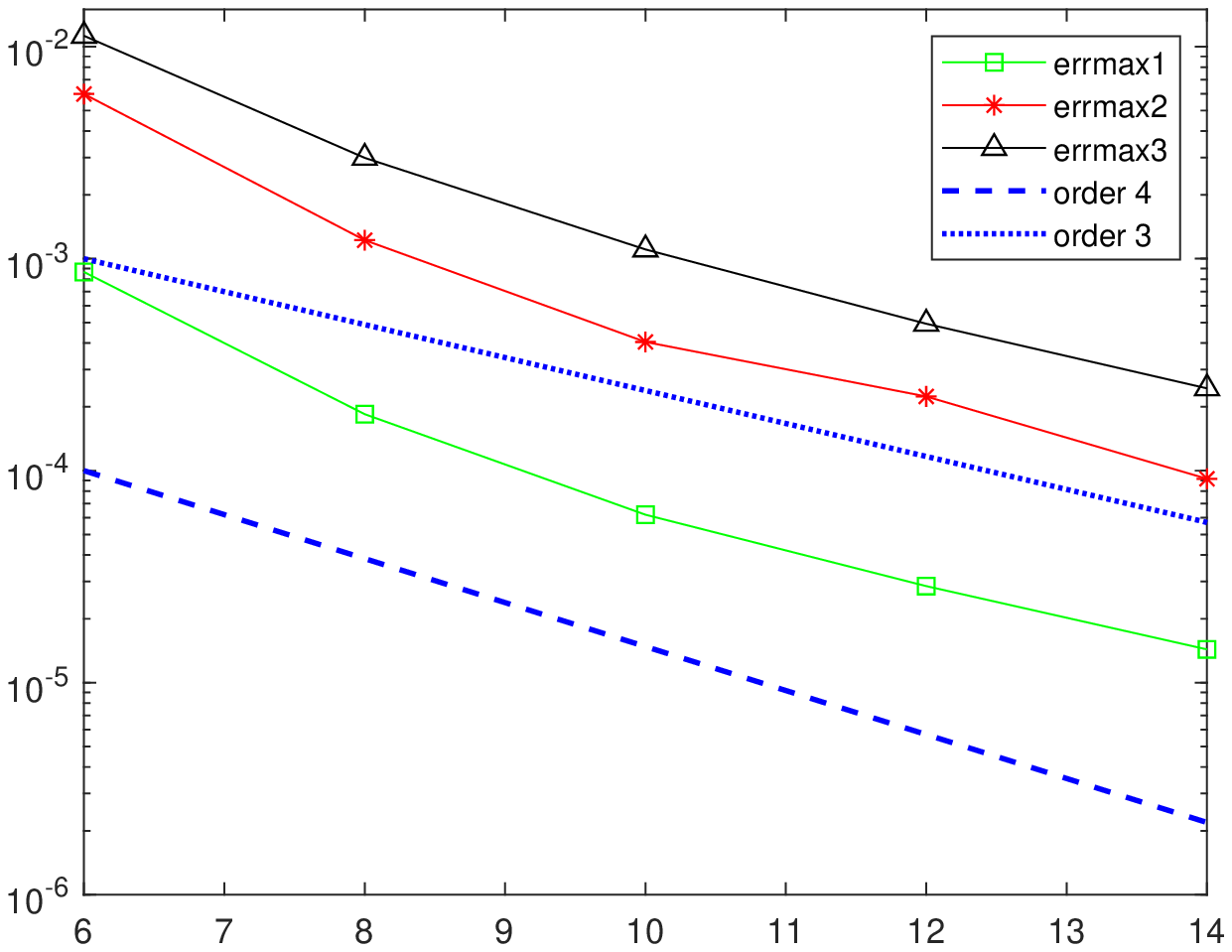}}
\subfigure[]{\includegraphics[width=6cm]{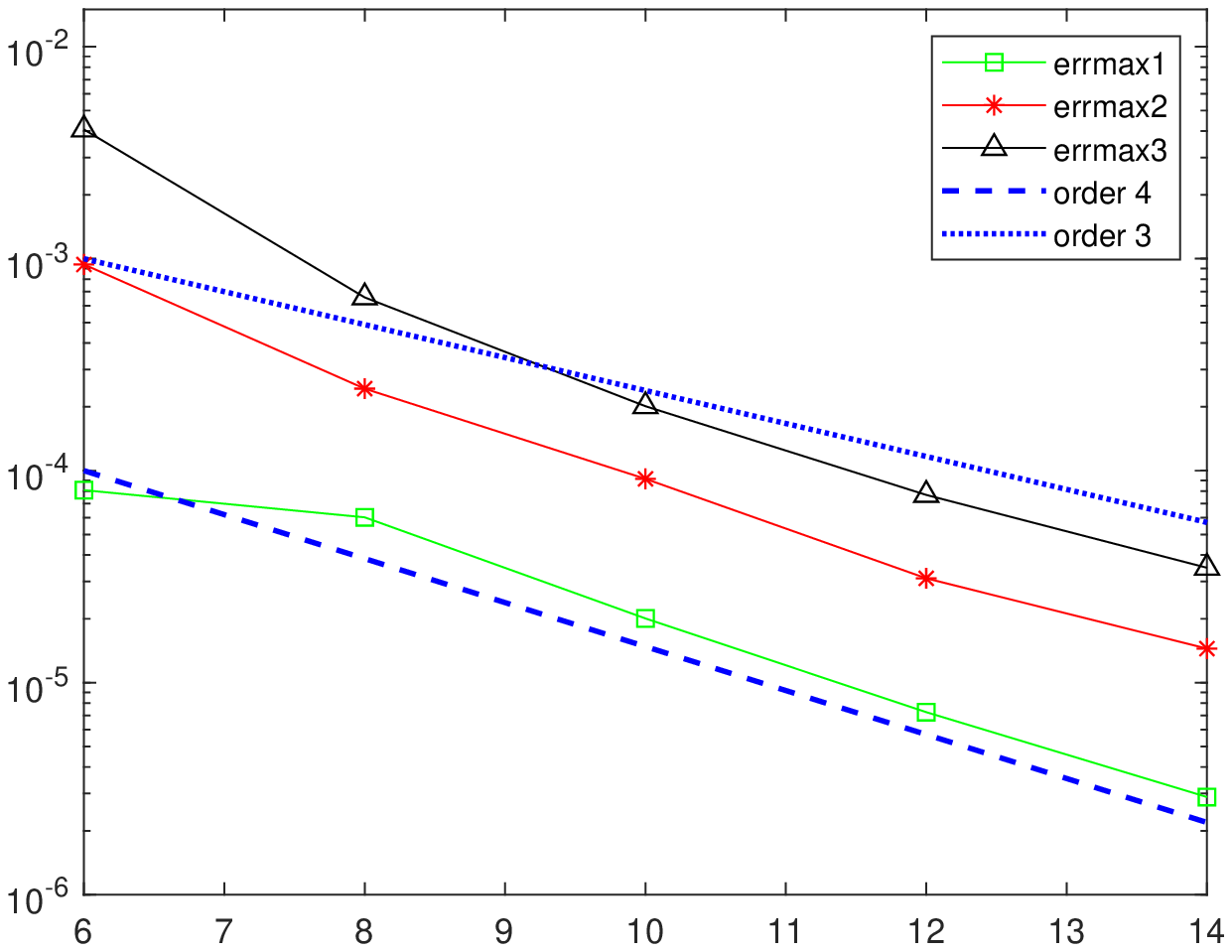}}
\subfigure[]{\includegraphics[width=6cm]{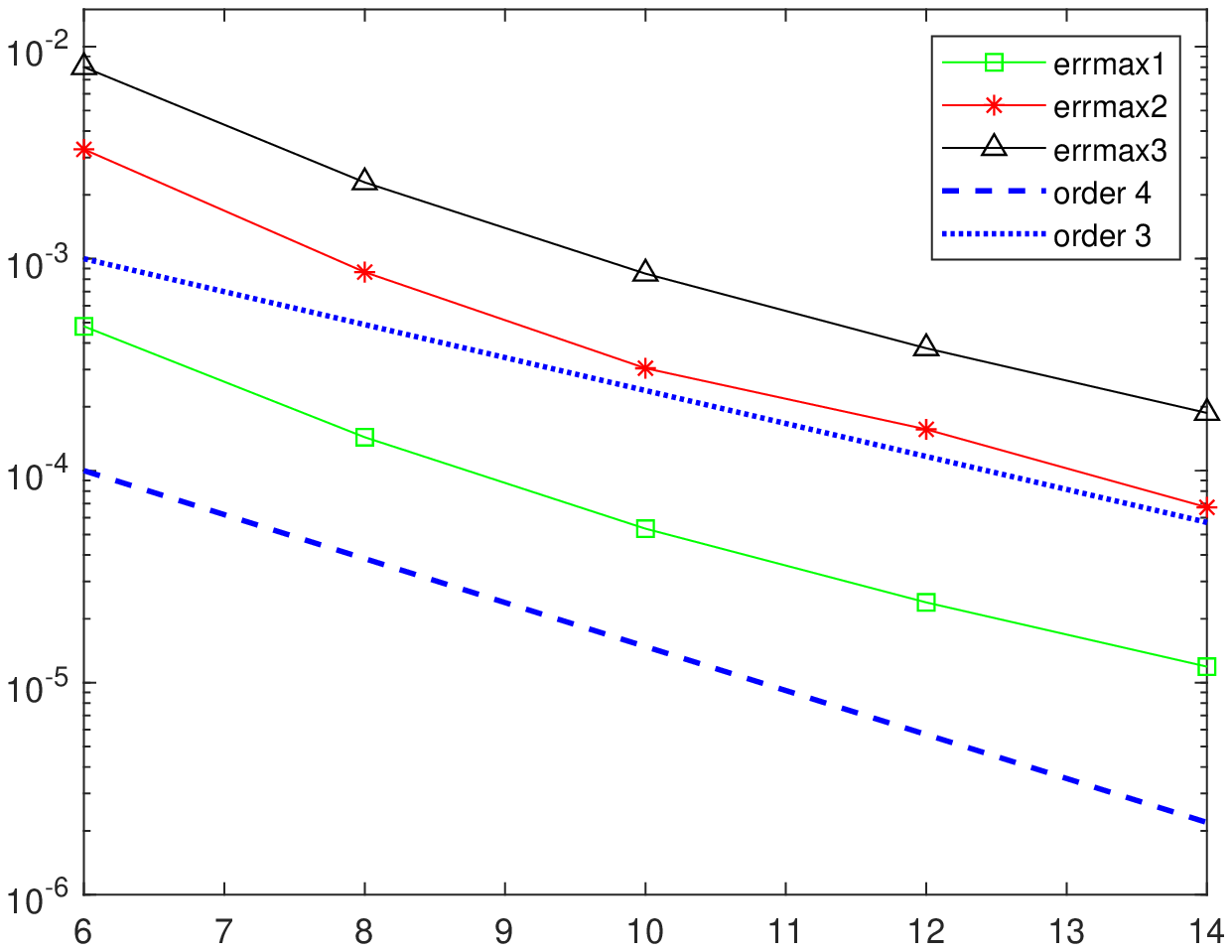}}
\subfigure[]{\includegraphics[width=6cm]{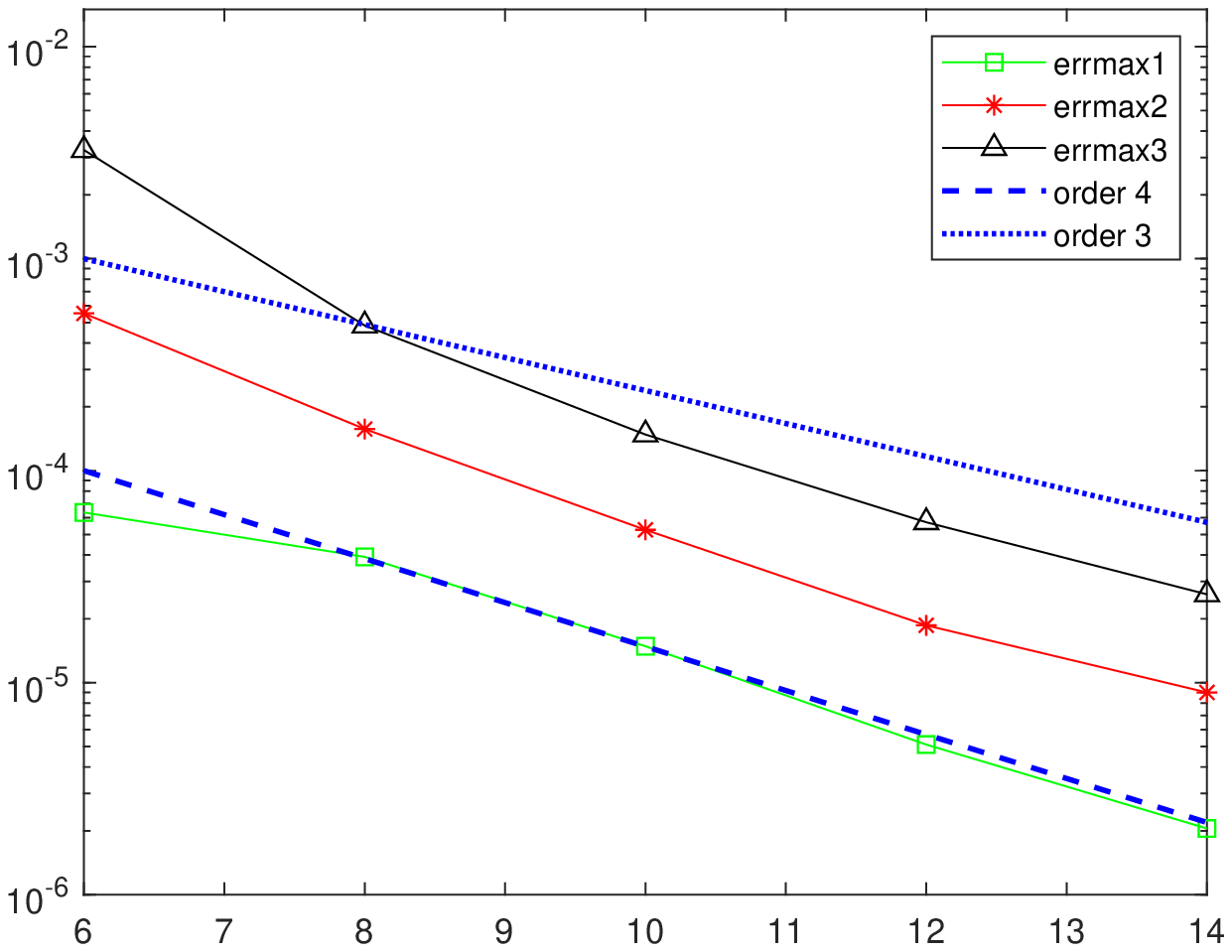}}

\caption{Example 4. The convergence behavior of the absolute cubature errors {\tt errmax1}, {\tt errmax2} and {\tt errmax3} for $d = 2, \ p = 2,$ (a),  $d = 2 , \ p = 3$ (b),
 $d = 3, \ p = 2$ (c) and  $d = 3, \ p = 3$ (d).}
\label{figEx4}
\end{figure}
As in the previous example, the integration domain $R_I$ is mapped to a quarter of the boundary of the considered section of hyperboloid whose height is $2.$
The matrix $A$ is again defined by the formula in (\ref{Adef}) but now the function $f_\s$ is assigned as follows,
$$f_\s(\t) = J(\t) \cos(k \Vert \X(\t)-\X(\s) \Vert_2)\,, $$
with $k = \pi/2$ and $J$ defined as in (\ref{Jdef}). Note that such function is $C^{1,1}$ also at $\t = \s.$
The so defined expression of (\ref{int}) is the real part of the weakly singular integral to be computed when the decomposition in (\ref{sum}) is applied for the Helmholtz kernel on the considered domain and the IgA--BEM collocation approach is adopted for the numerical solution.
The results for this example are shown in Figure \ref{figEx4}. From the figure we note that in this case increasing $p$ from $2$ to $3$ produced a better accuracy. The errors for the same value of $N$ are a  bit worse than those obtained in Example 3. This is due to the more oscillating nature of the function $f_\s.$ For a different approach to be applied in the nearly singular case with highly oscillating functions see  for instance \cite{Occorsio18}.

\section{Conclusions}
In this paper cubature rules for weakly singular double integrals containing an explicit B-spline factor are presented. The key ideas for these formulas are the extension   of a derivative free spline quasi-interpolation scheme and of an algorithm for spline product to the bivariate setting. Numerical results, also of interest in the IgA-BEM setting, confirm good performances of the proposed rules.

\section*{Acknowledgements}
The authors are all members of Gruppo Nazionale per il Calcolo Scientifico (GNCS) of the Istituto Nazionale di Alta Matematica (INdAM). The support of GNCS  through ``Progetti di ricerca 2019'' program is gratefully acknowledged. The first author is also thankful to the INdAM-GNCS funding ``Finanziamento Giovani Ricercatori 2020".
%
%

\end{document}